\documentclass[12pt,oneside]{amsart}
\usepackage{amsmath,amssymb,cite,mathrsfs,tikz-cd}
\usepackage[all]{xy}
\usepackage{typearea} 
\usepackage{hyperref} 
\usepackage{color} 
\usepackage[author={Max Schlepzig}]{pdfcomment}

\newcounter{maincounter}
\numberwithin{maincounter}{section}
\numberwithin{equation}{section}

\newtheorem{lemma}[maincounter]{Lemma}
\newtheorem{proposition}[maincounter]{Proposition}
\newtheorem{corollary}[maincounter]{Corollary}

\newtheorem{theorem}[maincounter]{Theorem}
\newtheorem{conjecture}[maincounter]{Conjecture}

\def\NN{\mathbb{N}}
\def\RR{\mathbb{R}}

\def\ZZ{\mathbb{Z}}
\def\PP{\mathbb{P}}
\def\QQ{\mathbb{Q}}

\newcommand{\cal}{\mathcal}

\newcommand{\cA}{\cal{A}}
\newcommand{\cP}{\cal{P}}

\newcommand{\cB}{\cal{B}}

\newcommand{\stab}[1]{{\mathrm{Stab}(#1)}}

\newcommand{\IP}{{\PP}}

\newcommand{\IR}{{\RR}}

\newcommand{\IZ}{{\ZZ}}
\newcommand{\IN}{{\NN}}

\newcommand{\IQ}{{\QQ}}

\newcommand{\cC}{{\mathcal C}}

\newcommand{\cJ}{{\mathcal J}}

\newcommand{\cO}{{\mathcal O}}

\newcommand{\cX}{{\mathcal{X}}}

\newcommand{\cZ}{{\mathcal Z}}

\newcommand{\ssm}{\setminus}

\renewcommand{\subset}{\subseteq} 

\newcommand{\jac}[1]{\mathrm{Jac}({#1})}
\newcommand{\diffS}{\mathrm{Diff}}

\newcommand{\tr}[2]{{\mathrm{Tr}}_{{#1}}({#2})}

\newcommand{\rk}[1]{\mathrm{Rk}({#1})}

\newcommand{\ntconst}{c_\mathrm{NT}}

\begin{document}
\title[Uniform Bound Pencil of Curves]{Uniform Bound for the Number of Rational Points on a Pencil
  of Curves}
\author{Vesselin Dimitrov}
\author{Ziyang Gao}
\author{Philipp Habegger}

\address{Department of Pure Mathematics and Mathematical Statistics, Centre for Mathematical Sciences, Wilberforce Road, Cambridge CB3 0WA, UK}
\email{vesselin.dimitrov@gmail.com}
\address{CNRS, IMJ-PRG, 4 place de Jussieu, 75005 Paris, France}
\email{ziyang.gao@imj-prg.fr}
\address{Department of Mathematics and Computer Science, University of Basel, Spiegelgasse 1, 4051 Basel, Switzerland}
\email{philipp.habegger@unibas.ch}

\subjclass[2000]{11G30, 11G50, 14G05, 14G25}

\maketitle
\begin{abstract}
Consider a one-parameter family of smooth, irreducible, projective curves of genus
$g\ge 2$ defined over a number field.
Each fiber contains at most finitely many rational points by
the Mordell Conjecture, a theorem of 
Faltings. We show that the number of rational points is bounded only
in terms of the family and the Mordell--Weil rank of the fiber's
Jacobian. Our proof uses Vojta's approach to the Mordell Conjecture
furnished with a height inequality due to the second- and third-named
authors. In addition we obtain uniform bounds for the number of torsion points in the
Jacobian that lie in each fiber of the family. 
\end{abstract}
\tableofcontents

\section{Introduction}

Let $k$ be a number field and $C$ a smooth,  irreducible,
projective curve of genus $g\ge 2$ defined over $k$. A fundamental
theorem of Faltings states that $C(k)$ is finite. Vojta
\cite{Vojta:siegelcompact} gave a different proof which was based
around an inequality of heights drawing on ideas from diophantine
approximations. No effective height upper bound for points in $C(k)$
is known. However, several authors including Bombieri
\cite{Bombieri:Mordell}, de Diego \cite{deDiego:97}, and R\'emond
\cite{Remond:Decompte} refined Vojta's approach to obtain estimates for
the cardinality of $C(k)$. In $\mathsection$\ref{sec:vojta} we review how to apply
Vojta's Method to curves.

Suppose $C$ has a $k$-rational point. We use it  to embed $C$ into
its Jacobian $\jac{C}$. Recall that by the Mordell--Weil Theorem
$\jac{C}(k)$ is a finite generated abelian group; we let $\rk{G}$
denote the  rank of any finitely generated abelian group $G$. We consider a suitable N\'eron--Tate height on
$\jac{C}$. By \cite[Theorem~2]{Bombieri:Mordell}, the number of points
in $C(k)$ of sufficiently large N\'eron--Tate height is bounded only
in terms of $\rk{\jac{C}(k)}$. The height cutoff depends on $C$ as
described by de Diego \cite{deDiego:97} . For a fixed curve, the set
of points in $C(k)$ whose height is below the cutoff is finite by
Northcott's Theorem. However, merely referring to Northcott's Theorem
induces a dependency on the curve. R\'emond \cite{Remond:Decompte}
gave a refined bound for the number of points of bounded height. His
approach relied on a lower bound for the N\'eron--Tate height such as
developed by David--Philippon \cite{DPvarabII}. But the cardinality
bound for $C(k)$ still depends on the Faltings height of the Jacobian
variety $\jac{C}$.
David and Philippon \cite{DaPh:07} considered a
different situation where $C$ is embedded in the power $E^g$ of an
elliptic curve $E$ defined over $k$. In this setting they proved a
height lower bound on $E^g$ with a correct dependency on the Faltings
height of $E$. As a consequence, David and Philippon
proved~\cite[Th\'eor\`eme~1.13]{DaPh:07} that if $C$ lies in $E^g$ and is not the translate of an
algebraic subgroup of $E^g$, then $\# C(k) \le c^{1+\rk{E(k)}}$
where $c>1$ is explicit and depends only on $g$ and suitable notion of degree of $C$.
Earlier, Silverman~\cite{Silverman:twists} produced a cardinality
estimate of the same quality when $C$ runs over twists of a fixed
smooth, projective curve of genus at least $2$. Here the rank of the
group of $k$-points of the
Jacobian of $C$ appears in the exponent. 
David, Nakamaye, and Philippon \cite[Th\'eor\`eme~1.1]{DaNaPh:07}
prove the existence of
a non-constant family of curves over $\IQ$ such that the number of $k$-rational
points on each fiber is bounded from above in terms of $g$ and $k$.

The purpose of this paper is to obtain cardinality estimates  of the
same quality as in the results of David--Philippon and Silverman
in a one-parameter family of curves. 
Our main tool is a
height lower bound of the second- and third-named authors
\cite[Theorem~1.4]{GaoHab} that replaces the result of David and
Philippon.

The following conjecture  is often attributed to
Mazur. He stated it as a question in a slightly less precise form  \cite[page~223]{Mazur:00} on bounding the number of rational points on curves.

\begin{conjecture}
  \label{conj:mazur} 
  Let $g\ge 2$ be an integer and $k$ a number field, there exists a
  constant $c=c(g,k)$ with the following property. If $C$ is a smooth
  projective curve of genus $g$ and defined over $k$, then the
  cardinality satisfies
  \begin{equation*}
    \# C(k) \le c^{1+\rk{\jac{C}(k)}}. 
  \end{equation*}
\end{conjecture}

Caporaso, Harris, and Mazur proposed the stronger Uniformity
Conjecture \cite{CaHaMa} where the upper bound is independent of
$\rk{\jac{C}(k)}$.



Let $\overline k$ be an algebraic closure of $k$. 
Prior to Conjecture~\ref{conj:mazur}, Mazur \cite[first~paragraph~on~page~234]{mazur1986arithmetic} asked a stronger but less formal question about the cardinality of the intersection of $C(\overline{k})$ with a finite rank subgroup of $\mathrm{Jac}(C)(\overline{k})$. Our main result is an attempt to answer this stronger question for any one-parameter family of
curves. 
  Let $S$ be a smooth,  geometrically irreducible
  curve defined over $k$. 
We suppose $S$ is embedded in some projective space and that
$h \colon S(\overline k)\rightarrow\IR$ is the pull-back of   the
absolute logarithmic Weil height.

\begin{theorem}
  \label{thm:main}
    Let $\cC$ be a geometrically irreducible, quasi-projective variety defined over $k$ together with a
  smooth morphism
 $\pi \colon \cC\rightarrow S$ with all fibers being
    smooth, geometrically irreducible,  projective curves of genus $g\ge 2$. Then 
   there exists a constant $c\ge 1$ depending on $\cC, \pi,$ and the choice
   of embedding of $S$ into projective space with the following
  property. Let $s\in S(\overline k)$ be such that $h(s)\ge c$ and
  suppose $\cC_s = \pi^{-1}(s)$ is
  embedded in its Jacobian $\jac{\cC_s}$ via the Abel--Jacobi map
  based at a $\overline k$-point of $\cC_s$. If $\Gamma$ is a finite rank subgroup of
  $\jac{\cC_s}(\overline k)$ with rank $\rk\Gamma$, then $\# \cC_s(\overline k)\cap \Gamma
  \le c^{1+\rk{\Gamma}}$. 
\end{theorem}

Note that in Theorem~\ref{thm:main}, the constant $c$ does not depend
on the choice of the $\overline k$-point via which the Abel-Jacobi
embedding is made. It is also independent of $\Gamma$.

Let $L$ be a finite extension of $k$, and let us now turn to
$L$-rational points on $\cC_s$. By Northcott's Theorem  the number of
$s\in S(L)$  with $h(s)<c$ is bounded from above in terms of
$c,[L:k],S,$ and the choice of projective embedding of $S$. The next
corollary follows directly  from Theorem~\ref{thm:main} applied to
$\Gamma = \jac{\cC_s}(L)$ combined with Faltings's Theorem for the
finitely many $s\in S(L)$ with $h(s)<c$. 
\begin{corollary}
  \label{cor:main}
  Let $\cC$ be as in Theorem~\ref{thm:main}, and let $L$ be a finite extension of $k$. 
   There exists a constant $c \ge 1$ depending on $\cC,\pi,$
 the choice
   of embedding of $S$ into projective space,  and $[L:k]$ with the following
  property. Let $s\in S(L)$, then $\# \cC_s(L)
  \le c^{1+\rk{\jac{\cC_s}(L)}}$.   
\end{corollary}


One can go a step further by applying R\'emond's
completely explicit \cite[Th\'eor\`eme~1.2]{Remond:Decompte},
\textit{cf.} \cite[page~643]{DPvarabII}, to handle the case $h(s)<c$.
Using the notation of \cite{DPvarabII} one can show that $h_0(\jac{\cC_s})$
is bounded from above linearly in terms of $[L:k] \max\{1,h(s)\}$. 
So we may choose $c$ in Corollary \ref{cor:main} to depend polynomially on $[L:k]$.


Our approach to bounding the number of rational points is ultimately
based on Vojta's Method. 
Recently, Alpoge \cite{AlpogeRatPt}, based on this method, proved the following result: for genus $2$ curves with a marked Weierstrass point, the average number of rational points is bounded. 
We recall some other approaches, one going
back to work of Chabauty \cite{Chabauty}, under an additional
hypothesis on the rank of Mordell--Weil group. Let $C$ again be a
smooth, geometrically irreducible, projective curve $C$ defined over $k$ and of
genus $g\ge 2$. The Chabauty--Coleman approach
\cite{Coleman:effCha,Stoll:Uniform,KatzRabinoffZB} leads to strong
bounds for $\#C(k)$ if the rank of $\jac{C}(k)$ is small in terms of
$g$. For example, if $C$ is hyperelliptic and
$\jac{C}(k)$ has rank at most $g-3$, then Stoll
\cite{Stoll:Uniform} showed that the cardinality of $C(K)$
is bounded solely in terms of $K$ and $g$.
Stoll's approach inspired the work of Katz--Rabinoff--Zureick-Brown
\cite{KatzRabinoffZB} who proved that the cardinality is bounded  only
in terms of  $g$
and $[K:\IQ]$ and without  the hyperelliptic
hypothesis. 
An approach based on connections to
unlikely intersections was investigated by Checcoli, Veneziano, and
Viada \cite{CVV:17} again under the assumption that the rank of
Mordell--Weil group of an ambient abelian variety is sufficiently
small in terms of the dimension. Our result does not stipulate a
restriction on $\rk{\jac{C}(k)}$, but it is confined to a one-parameter
family of curves. Recall that the coarse moduli space of all genus
$g\ge 2$ curves has dimension $3g-3\ge 3$.

Recall that Raynaud's Theorem, the
Manin--Mumford Conjecture, states that the image of $C$ under an
Abel--Jacobi map $C\rightarrow \jac{C}$  meets at most finitely many
points in the full torsion subgroup
$\jac{\cC_s}_{\mathrm{tor}}$ of $\jac{\cC_s}(\overline k)$.
Theorem~\ref{thm:main} applied to $\Gamma =\jac{\cC_s}_{\mathrm{tor}}$  leads to the following uniform bound in a
one-parameter family.
\begin{corollary}
  Let $\cC$ and $c$ be as in Theorem~\ref{thm:main}. Let $s\in
  S(\overline k)$ be such that $h(s)\ge c$ and suppose $\cC_s$ is
  embedded in its Jacobian $\jac{\cC_s}$ via the Abel-Jacobi map based
  at a $\overline k$-point of $\cC_s$. Then $\#\cC_s(\overline k)\cap
  \jac{\cC_s}_{\mathrm{tor}}\le c$. 
\end{corollary}


As in Corollary \ref{cor:main} we can apply Northcott's Theorem if we restrict $s$ to $S(L)$ with $L$ a finite extension of $k$.
The following corollary is a direct consequence of Theorem \ref{thm:main} together
with David and Philippon's \cite[Th\'eor\`eme~1.2]{DPvarabII}, see also R\'emond 
\cite[Th\'eor\`eme~1.2]{Remond:Decompte},  that handles 
the finitely many points $s\in S(L)$ with $h(s)<c$ and gives a
cardinality bound  independent of the base point.
\begin{corollary}
  \label{cor:torsion}
  Let $\cC$ be as in Theorem~\ref{thm:main}, and let $L$ be a finite extension of $k$. 
   There exists a constant $c\ge 1$ depending on $\cC,\pi,$ 
 the choice
   of embedding of $S$ into projective space, and $[L:k]$ with the following
   property. Let $s\in S(L)$ and
suppose $\cC_s$ is
  embedded in its Jacobian $\jac{\cC_s}$ via the Abel--Jacobi map
  based at a $\overline k$-point of $\cC_s$.
     Then 
  $\cC_s(\overline k)$ contains  at most $c$ torsion points of $\jac{\cC_s}$. 
\end{corollary}
We can use David and Philippon's  \cite[Th\'eor\`eme~1.2]{DPvarabII}
to produce a constant $c$ that depends polynomially on $[L:k]$.
Indeed, in their Th\'eor\`eme 1.4 and their notation $q(\cC_s)$ is bounded
polynomially in $[L:k]$ if $h(s)<c$ when $c$ is fixed.

Recently, DeMarco--Krieger--Ye \cite{DeMarcoKriegerYeUniManinMumford}
proved a bound on the cardinality of torsion points on any genus $2$
curve that admits a degree-two map to an elliptic curve when the
Abel--Jacobi map is based at a Weierstrass point. Their result is not
confined to one-parameter families and their bound $c$ is furthermore
independent of $[L:k]$, or stronger, independent of $[L:\IQ]$.

In future work, we plan to develop the height bound \cite{GaoHab}
beyond base dimension one and to use the approach presented in this
paper to generalize Theorem~\ref{thm:main} accordingly.


\subsection*{Acknowledgements} The authors thank Gabriel Dill, Fabien Pazuki,
Patrice Philippon, and Ga\"el R\'emond for helpful comments on a
preliminary version of this text.  They also thank Christian Elsholtz
for comments and
for suggesting  the application involving the Hardy--Ramanujan Theorem. They also thank the referees for their careful reading and valuable comments. Vesselin Dimitrov would like to thank the NSF and the Giorgio and Elena Petronio Fellowship Fund II for financial support for this work. Both VD and ZG would like to thank the Institute for Advanced Study and the special year ``Locally Symmetric Spaces: Analytical and Topological Aspects'' for its hospitality during this work.

\section{A review of Vojta's Method}
\label{sec:vojta}

In this section we give an overview, at least up to the constants
involved, of bounding the cardinality in the Mordell Conjecture using
Vojta's approach. No new material is contained here. 
We rely on R\'emond's \cite{Remond:Vojtasup,Remond:Decompte}
quantitative results. Work of Pazuki
\cite{Pazuki:uniform} also
involves completely explicit constants.

Let $k$ be a number field with a fixed algebraic closure $\overline
k$ of $k$. Let $A$ be an abelian variety defined over $\overline k$ equipped
with a very ample and symmetrical line bundle.
We may suppose that
 an  immersion   attached to the line bundle
 realizes $A$ as a  projectively
normal subvariety of $\IP^n$. 

Let $h$ denote the absolute logarithmic Weil height on
$\IP^n(\overline k)$.
We write $h_1$ for an upper bound for the absolute logarithmic
projective height of bihomogeneous polynomials that
describe the addition morphism on $A$ as a subvariety of $\IP^n$, see
Section 5~\cite{Remond:Vojtasup} for a precise definition. 
Tate's Limit Process provides us with a
N\'eron--Tate height $\hat h \colon A(\overline k)\rightarrow\IR$. It is well-known
 there exists a constant $\ntconst$, depending
on the data introduced above, such
that $|h(P)- \hat h(P)|\le \ntconst$ for all $P\in A(\overline
k)$. For $P,Q\in A(\overline k)$  we set
$\langle P,Q\rangle = (\hat h(P+Q)-\hat h(P)-\hat h(Q))/2$
and often abbreviate $|P| = \hat h(P)^{1/2}$. The notation $|P|$ is
justified by the fact that it induces a norm after tensoring with the reals.

Any irreducible closed subvariety $X$ of $\IP^n$ has a well-defined
degree $\deg X$ and height $h(X)$; for the latter we refer to \cite{BGS}.

A \textit{coset} in an abelian variety is a translate of an abelian
subvariety.
A coset is called \textit{proper} if it is not equal to the ambient
abelian variety. 

For what follows let $C$ denote an irreducible curve contained in $A$. 
The following theorem is a special case of  R\'emond's
\cite[Th\'eor\`eme~1.2]{Remond:Vojtasup}.
It is a version of Vojta's Inequality.

\begin{theorem}
  \label{thm:vojta}
  There exists a constant $c=c(n)\ge 1$ depending only on $n$ with the
  following property. Let
    $c_1 = c\deg(C)^2$ and $c_2=c \deg(C)^6$
and suppose $C$ is not a coset in $A$.
If $P,Q\in C(\overline k)$ satisfy
\begin{equation*}
 \langle P,Q\rangle \ge \left(1-\frac{1}{c_1}\right) |P||Q|
\quad\text{and}\quad |Q| \ge c_2 
|P|
\end{equation*}
then 
\begin{equation*}
|P|^2 \le c \deg(C)^{20}
    \max\{1,h(C),h_1,\ntconst\}.
\end{equation*}
\end{theorem}

A second tool is the so-called Mumford equality. We  use a
quantitative version due to R\'emond. We write $\stab{C}$ for the
stabilizer of $C\subset A$, it is an algebraic subgroup of $A$. 

\begin{theorem}
\label{thm:mumford}  There exists a constant $c=c(n)\ge 1$ depending only on $n$ with the
  following property. 
Say  $P,Q\in C(\overline k)$ with $P-Q \not \in \stab{C}(\overline k)$. If 
  \begin{equation*}
 \langle P,Q\rangle \ge \left(1-\frac{1}{c \deg(C)^2}\right)  |P|
 |Q| \quad\text{and}\quad
 \bigl| |P| -  |Q|\bigr| \le \frac{1}{c \deg(C)}  |P|
  \end{equation*}
  then
  \begin{equation*}
|P|^2  \le c \deg(C)^3 \max\{1,h(C), h_1,\ntconst\}.
  \end{equation*}
\end{theorem}
\begin{proof}
  In R\'emond's \cite[Proposition~3.4]{Remond:Decompte}
  we take $c_1=c \deg(C)^2$ and $c_4=c \deg(C)$ 
where  $c$ is large enough in terms of $n$.
  The condition that $P-Q$ is not in the stabilizer of $C$ is
implies that $(P,Q)$ is isolated in the fiber of the
  subtraction morphism $A\times A\rightarrow A$ restricted to $C\times
  C$. 
\end{proof}

The precise exponent of $\deg C$ in the results above is irrelevant
for our main results. 

It is well-known how the inequalities of Vojta and Mumford combine to
yield the following result. For the reader's convenience we recall here
this classical argument.

\begin{corollary}
  \label{cor:countlargepts}
  There exists a constant $c=c(n,\deg C)\ge 1$ depending
  only on $n$ and $\deg C$ with the following property. 
  Suppose $\Gamma$ is a subgroup of $A(\overline k)$ of finite rank
  $\rho \ge 0$. 
  If $C$ is not a coset in $A$, then
  \begin{equation*}
   \# \left\{ P \in C(\overline k)\cap \Gamma : |P|^2 > c
   \max\{1,h(C),h_1,\ntconst\} \right\} \le c^{\rho}. 
  \end{equation*}
\end{corollary}
\begin{proof}
The hypothesis on $C$ implies that there exist $R,R'\in C(\overline
k)$ with $C-R \not = C-R'$. The stabilizer $\stab{C}$ lies in the
finite set $(C-R)\cap (C-R')$ which has cardinality at most
$\deg(C)^2$ by a suitable version of B\'ezout's Theorem.

Observe that both Theorems~\ref{thm:vojta} and \ref{thm:mumford}
hold with $c$ replaced by some larger value.  
We let $c$ denote the maximum of both constants $c$ from these two
theorems.
  
Let $P_1,P_2,\ldots,P_N \in C(\overline k)\cap \Gamma$ be pairwise
distinct points such that
\begin{equation}
  \label{eq:Piincreasing}
   c \deg(C)^{20} \max\{1,h(C),h_1,\ntconst \} < |P_1|^2 \le |P_2|^2
 \le \cdots.
  \end{equation}

For given $P_i$ there are at most   $\#\stab{C}\le \deg(C)^2$
different $P_j$ with $P_i-P_j \in \stab{C}(\overline k)$. By the
pigeonhole principle there are  $N'\ge N/\deg(C)^2$ members among
$P_i$ 
whose pairwise difference is $0$ or not in $\stab{C}(\overline k)$.
After thinning out our sequence and renumbering we may assume, in
addition to (\ref{eq:Piincreasing}) that 
\begin{equation}
  \label{eq:PiPjnotinstab}
  P_i-P_j\not\in
\stab{C}(\overline k)
\end{equation}
 for all $i,j\in \{1,\ldots,N'\}$ with $i\not=j$. 

The $\rho$-dimensional $\IR$-vector space $\Gamma\otimes_\IZ\IR$ is equipped with an inner product induced by 
the N\'eron--Tate pairing $\langle \cdot,\cdot\rangle$.
We also write $|\cdot|$ for the resulting norm on
$\Gamma\otimes_\IZ\IR$. 
This
allows us to do Euclidean geometry in $\Gamma\otimes_\IZ\IR$. It is no restriction to assume
$\rho\ge 1$. By R\'emond's \cite[Corollaire~6.1]{Remond:Decompte}, the
vector space can be covered by at most  $\lfloor
(1+({8c_1})^{1/2})^\rho\rfloor$ cones on which $\langle P,Q\rangle \ge
(1-1/c_1) |P| |Q|$ holds where  for
$c_1$ we pick the constant from Theorem~\ref{thm:vojta}. 
We use again the pigeonhole principle to thin out
$P_1,\ldots,P_{N'}$ and get new subsequence of $N'' \ge N'/
(1+({8c_1})^{1/2})^\rho$ pairwise distinct points that lie in a common cone.
Thus
\begin{equation}
  \label{eq:PiPjincone}
  \langle P_i,P_j\rangle \ge \left(1-\frac{1}{c_1}\right) |P_i||P_j|
\end{equation}
for all $i,j\in \{1,\ldots,N''\}$.

By
(\ref{eq:Piincreasing}) and
(\ref{eq:PiPjincone}) the hypothesis in Theorem~\ref{thm:vojta} cannot be
true. Thus
 we must have
\begin{equation*}
|P_1|\le  |P_i| \le c_2 |P_1|
\end{equation*}
for all $i\in \{1,\ldots,N''\}$.
We apply the pigeonhole principle a
final time. Let us assume $N''\ge 2$. Considering
the sequence $|P_1| \le \cdots \le |P_{N''}|$ there must exist a
consecutive pair $P=P_i$ and $Q=P_{i+1}$ with
\begin{equation*}
 \bigl||P| -|Q|\bigr| \le  \frac{1}{N''-1} (|P_{N''}|-|P_1|)\le \frac{c_2}{N''-1} |P_1|.
\end{equation*}

By (\ref{eq:PiPjnotinstab}) we know that $P-Q \not \in
\stab{C}(\overline k)$.
Furthermore,
by (\ref{eq:Piincreasing}) and (\ref{eq:PiPjincone}), the hypothesis
of Theorem~\ref{thm:mumford} cannot be met, thus
$c_2 / (N''-1) > 1 /(c\deg C)$ and
$N'' < c c_2 \deg(C) + 1$. 

The last paragraph implies $N'' \le \max\{2,c c_2 \deg(C)+1\}$. Recall
that $c$ and $c_2$ depend only on $n$ and $\deg(C)$. The
corollary follows from 
\begin{equation*}
  N \le \deg(C)^2 N' \le (1+({8 c_1})^{1/2})^\rho \deg(C)^2 N''.\qedhere
\end{equation*}
\end{proof}

In Corollary \ref{cor:countlargepts2}  below we subsum these important results in the relative setting of
an abelian scheme over a curve. 
Let $S$ be a smooth,  irreducible affine curve  defined
over $\overline k$ and let 
$\pi \colon \cA\rightarrow S$ be an abelian scheme. 
Let $\iota \colon \cA \rightarrow \IP^n \times S$ be an $S$-immersion such that the restriction of $\iota^*\cO_S(1)$ to $\cA_s = \pi^{-1}(s)$ is symmetric and 
$\cA_s\rightarrow \IP^n$ is projectively normal for all $s\in S(\overline k)$. We write $\hat h \colon \cA(\overline k)\rightarrow [0,\infty)$ for the fiberwise N\'eron--Tate height. 

 Moreover, we may assume that $S$ is contained in
some projective space. Then we  restrict the Weil height from
projective space and obtain a height function  $h \colon S(\overline
k)\rightarrow \IR$.

If $\cC$ is an irreducible closed subvariety of $\cA$ and $s\in
S(\overline k)$, then we write
$\cC_s$ for $\pi|_{\cC}^{-1}(s)$. See also de Diego 
\cite[Th\'eor\`eme~2]{deDiego:97}.

We may reformulate Corollary~\ref{cor:countlargepts} as follows. 
\begin{corollary}
    \label{cor:countlargepts2}
  Let $\cC \subset \cA$ be an irreducible closed subvariety that
  dominates $S$ and such that the generic fiber  of $\cC\rightarrow S$
  is a geometrically irreducible curve. 
  There exists a constant $c=c(n,\pi,h,\cC)\ge 1$ with the following property.
  Suppose $s\in S(\overline k)$ and that 
 $\Gamma$ is a subgroup of $\cA_s(\overline k)$ of finite rank
  $\rho \ge 0$. If $C$ is an irreducible component of $\cC_s$
  that is not a coset in $\cA_s$, then 
  \begin{equation*}
   \# \left\{ P \in C(\overline k)\cap \Gamma : \hat h(P) > c
   \max\{1,h(s)\} \right\} \le c^{\rho}. 
  \end{equation*}
\end{corollary}
\begin{proof}
  First we observe that each fiber $\cC_s$ is equidimensional of
  dimension $1$ and that the degree on any irreducible component
  is bounded from above  uniformly in $s$.
  
  To apply Corollary~\ref{cor:countlargepts} it is enough to note that
  $h(\cC_s), h_1,$ and $\ntconst$, all functions of $s$, are 
  bounded from above linearly in $\max\{1,h(s)\}$.

  We may regard $\cC$ as a subvariety of $\IP^n\times S$ and $S$ as
  inside some projective space. The fiber
  $\cC_s$ is the intersection of $\cC$ with $\IP^n\times\{s\}$. So the
  desired bound for $h(\cC_s)$  follows  from an appropriate version
  of the Arithmetic B\'ezout Theorem such as Th\'eor\`eme
  3~\cite{HauteursAlt3}.

  To bound $h_1$ note that polynomials defining addition exist on the
  generic fiber of $\cA\rightarrow S$ taken as a subvariety of $\IP^n$
  over the function field $\overline k(S)$. Specializing outside a
  finite set of $S(\overline k)$ gives desired polynomials on each
  fiber $\cA_s$. The finitely many exceptions can be compensated by
  increasing the linear factor in front of $\max\{1,h(s)\}$.

 There are several approaches to bounding $\ntconst$ from above. For
 example, we can apply the Theorem of
  Silverman--Tate~\cite[Theorem~A]{Silverman}. Indeed, by
  compactifying and then desingularizing using Hironaka's Theorem we can
  find a smooth, projective model for $\cA$ over a smooth, projective base as in the beginning of Section
  3~\cite{Silverman}. The structural morphism on this model is automatically flat as
  $\dim S = 1$.\footnote{In our special case, \textit{i.e.} if $\dim
    S = 1$, is likely that a suitable
  version of the Silverman--Tate Theorem holds where 
  Weil divisors are replaced by Cartier divisors and where the model
  for $\cA$ is merely normal. As desingularizing a curve is
  elementary,   this approach would  make invoking Hironaka's Theorem unnecessary.}
  Bounds for  $\ntconst$ were  studied by Manin and
  Zarhin~\cite{ZarhinManin}. 
 David and Philippon later gave an explicit estimate in Proposition
 3.9~\cite{DPvarabII} using Mumford coordinates. 
\end{proof}

When $\Gamma = \jac{C}(k)$ and the $C$ is embedded into $\jac{C}$ via a $k$-point, Alpoge \cite{AlpogeRatPt} has more explicit bounds for $c$.

\section{A preliminary lemma}
\label{sec:avoid}

Let $\overline k$ be an algebraically closed field and suppose $A$ is an
abelian variety  defined over $\overline k$. 
We define a homomorphism
$\diffS \colon A^3\rightarrow A^2$  of algebraic groups defined by
\begin{equation*}
  \diffS(P_0,P_1,P_2) = (P_1-P_0,P_2-P_0)
\end{equation*}
for all $(P_0,P_1,P_2)\in A^{3}{(\overline k)}$. This morphism
is sometimes called the Faltings--Zhang map and generalizes the
difference morphism studied by Bogomolov in
\cite[$\mathsection$2.8]{Bogomolov}.

Suppose $C\subset A$ is an irreducible closed curve; we will impose
additional conditions on $C$ below. 
The image
$\diffS(C^3)$ plays an important role in the proof of our main theorem.
Later on $A$ will arise as a fiber of an abelian
scheme and $C$ appears in a fiber of a family of curves. We will use a height inequality which holds on $\diffS(C^3)$
apart from a finite collection of naturally defined subvarieties.
In the current section we analyze this collection and how it relates to $C$.

In this section $H$ denotes an
abelian subvariety of $A$ that we should think of as  arising from the constant
part of the abelian scheme. 


\begin{lemma}\label{LemmaAvoid1}
Suppose $C$ is not a coset in $A$.
Let $B$ be an abelian subvariety of $A^2$. Suppose $Z$ is an irreducible closed subvariety of $H^2$ with
\[
(C-P)^2 = Z+B+(Q_1,Q_2)
\]
for some $P \in C(\overline k)$ and $(Q_1,Q_2) \in A^2(\overline k)$.
Then $C$ is contained in a translate of $H$ in $A$.
\end{lemma}
\begin{proof}
Denote $q \colon A^2 \rightarrow A$ the projection to the first factor. The assumption of the lemma implies $q(Z) + q(B) + Q_1 = C-P$. So the stabilizer of $C$ contains the abelian subvariety $q(B)$ of $A$. If $q(B) \not= 0$, then by dimension reasons we have that $C$ is a coset in $A$, contradicting our assumption in $C$. Hence $q(B) = 0$ and from above we find $C-P = q(Z) + Q_1 \subset q(H^2) + Q_1 = H+Q_1$. So $C \subset H+Q_1+P$ is contained in a translate of $H$ in $A$, as desired. 
\end{proof}

\begin{lemma}
  \label{lem:avoid}
Suppose $C$ is not a coset in $A$.
Let $B$ be an abelian subvariety of $A^2$ and $(Q_1,Q_2)\in
A^2(\overline k)$ with $(Q_1,Q_2)+B\subset
\diffS(C^{3})$. Suppose  
$(P_1,P_2)\in C(\overline k)^2$ such that
\begin{equation}
  \label{eq:difffiber}
 \diffS(C\times \{(P_1,P_2)\})\subset  B+(Q_1,Q_2). 
\end{equation}
Then $P_1-P_2=Q_1-Q_2$.
\end{lemma}
\begin{proof}
The condition (\ref{eq:difffiber}) implies 
\begin{equation}
    \label{eq:difffiber2}
(P_1,P_2)-(P_0,P_0)= (P_1-P_0,P_2-P_0) \in B(\overline k)+(Q_1,Q_2)
  \end{equation}  
for all    $P_0 \in
C(\overline k)$.
Let $\Delta:A\rightarrow A^2$ denote the diagonal embedding. Then
$\Delta(C) \subset (P_1-Q_1,P_2-Q_2)  + B$
and so $C$ lies in $\Delta^{-1}(((P_1-Q_1,P_2-Q_2)+B)\cap\Delta(A))$
which is
 a finite union of cosets in $A$ of  dimension $\dim B\cap \Delta(A)$.
But $C$ is no coset, thus
\begin{equation}
  \label{eq:dimBlb}
\dim B\ge \dim B\cap \Delta(A) \ge 2. 
\end{equation}

Let us assume for the moment that $\dim B=2$. 
From (\ref{eq:dimBlb})  we conclude 
$B \subset \Delta(A)$ since $B$ is irreducible. 
Fix a $\overline k$-point $P_0\in C(\overline k)$. By
(\ref{eq:difffiber2}) we have
$(P_1-P_0,P_2-P_0) \in  (Q_1,Q_2)+\Delta(A)(\overline
k)$
and so $(P_1-P_0,P_2-P_0) = (Q_1+Q,Q_2+Q)$ for some $Q\in A(\overline
k)$.
We cancel $Q$ by subtracting and find 
 $P_1-P_2=Q_1-Q_2$, as desired.

Now let us assume $\dim B\ge 3$. We will arrive at a contradiction.
Let $X = \diffS(C^3)$. This is an irreducible variety and $\dim X\le
3$. Then $(Q_1,Q_2)+B$, being irreducible and inside $X$ by hypothesis, must equal $X$ and $\dim B = \dim X = 3$. As above
 $q \colon A^2\rightarrow A$ denotes the projection onto the first factor. 
Any fiber of $q|_{B}$ can be identified with a
coset in the second factor of $A^2$. Moreover for any $P'_0\in C(\overline k), P'_1\in C(\overline k)$, we have that
$\{(P'_1-P'_0-Q_1)\}\times (C-P'_0-Q_2) \subset B$ is a curve contained in a
fiber of $q$. As $C$ is not a coset in $A$, some
fiber of $q|_B$ has dimension at least $2$. But then all fibers of
$q|_B$ have dimension at least $2$ and hence 
$\dim q(B)\le \dim B - 2= 1$. On the other hand $Q_1+q(B) = q(X)=C-C$ where
$(C-C)(\overline k)=\{P'_1-P'_0 : P'_0\in C(\overline k),P'_1\in
C(\overline k)\}$. In  particular, $\dim (C-C) = 1$ and therefore
$C- P= C-C$ for all $P\in C(\overline k)$.    Thus the stabilizer of
$C$ contains $C-C$ and in particular, $C$
is a coset in $A$. This contradicts the hypothesis.
  \end{proof}

Now we end this section with the following corollary. Note that $(C-P)^2 = \diffS(\{P\} \times C^2)$ for any $P \in C(\overline k)$.

\begin{corollary}\label{LemmaAvoid}
Suppose that $C$ is  not contained in any proper coset in $A$ and $C\not=A$.
  Let 
$\psi \colon A\rightarrow H$ be a homomorphism such that $\psi|_H \colon H\rightarrow
H$ is an isogeny. Let
$\psi\times\psi \colon A^2\rightarrow H^2$ be the square of $\psi$ and let
$B$ be an abelian subvariety of $A^2$ contained  in the kernel of    $\psi\times\psi$.
Suppose  $Z$ is an irreducible closed subvariety of $H^2$ and $Q_1,Q_2\in A(\overline k)$ such that 
$Z+B+(Q_1,Q_2) \subset \diffS(C^3)$.
If $H\subsetneq A$, then 
 $(C-P)^2 \not\subset Z+B+(Q_1,Q_2)$
for all $P\in C(\overline k)$.
\end{corollary}
\begin{proof}
By hypothesis $C$ is not a coset in $A$ and not
contained in a translate of
  $H$ in $A$.
  Let us assume $H\not=A$ and that
  there exists $P\in C(\overline k)$ with
$(C-P)^2 \subset Z+B+(Q_1,Q_2)$. To prove the corollary we will derive
  a contradiction.

For dimension reasons, $Z+B+(Q_1,Q_2)$ is either $(C-P)^2$ or
$\diffS(C^3)$. The first case is impossible by
Lemma~\ref{LemmaAvoid1}. So
\begin{equation}
  \label{eq:diffCZB}
 \diffS(C^3)=Z+B+(Q_1,Q_2).
\end{equation}
  As $C$ is not contained in a translate of $H$ there exist
  $P_1,P_2\in C(\overline k)$ such that $P_1-P_2 \not\in H(\overline k) + (Q_1-Q_2)$. 
In particular, 
\begin{equation}
  \label{eq:negi}
  \diffS(C\times \{(P_1,P_2)\})\subset Z+B + (Q_1,Q_2).
\end{equation}

We claim that $Z$ is not a coset in
$H^2$. Indeed, otherwise $Z=(Q_1',Q_2') + B'$ for an abelian subvariety
$B'\subset H^2$ and $Q_1',Q_2'\in H(\overline k)$. 
 In this case, we get a contradiction from
(\ref{eq:negi}) and Lemma~\ref{lem:avoid} applied to the abelian subvariety $B'+B$ of $A^2$
 and using $P_1-P_2\not\in H(\overline k)+(Q_1 - Q_2)$.
In particular, $Z\subsetneq H^2$ and $\dim H\ge 1$. 

Let us suppose  $\dim H = 1$. 
Recall that $\psi(C)\subset H$ and $\psi(C)$ cannot be a point, so $\psi(C)=H$. 
Recall also $B \subset \ker \psi\times\psi$  so (\ref{eq:negi}) implies that
$\psi\times\psi(\Delta(C))$ is contained in a translate of 
$-\psi\times\psi(Z)$.
So $\Delta(H)=\Delta(\psi(C))=\psi\times\psi(\Delta(C))$ lies  in a
translate of $-\psi\times\psi(Z)$. 
The last paragraph  implies $\dim\psi\times\psi(Z)\le \dim Z\le  \dim H^2-1=
1$.  Because $\dim\Delta(H)=\dim H= 1$ we see that  $\psi\times\psi(Z)$   is a translate of
the diagonal  $\Delta(H)$. In particular, 
 $\psi\times\psi(Z)$ is a coset in 
 $H^2$.
But $Z\subset H^2$ and $\psi\times\psi|_{H^2}\colon H^2\rightarrow H^2$ is an isogeny
by hypothesis.
Therefore, $Z$ is a coset in $H^2$, something we excluded above. 

 Hence $\dim H\ge 2$ and thus  $\dim A \ge 3$ as $H\subsetneq A$.  In particular, $C-C$ is not an abelian surface, as otherwise $C$ would be contained in a proper coset in $A$. So $\dim \mathrm{Stab}(C-C) = 0$ or $1$.

Suppose $\dim \mathrm{Stab}(C-C) = 0$. 
Recall that $q \colon A^2 \rightarrow A$ is the projection to the
first factor. Then $C - C=q(\diffS(C^3))=q(Z) + q(B) + Q_1$ by (\ref{eq:diffCZB}). Since
$C-C$ has trivial stabilizer, we have $q(B) = 0$ and thus $C-C =
q(Z) + Q_1 \subset q(H^2) + Q_1 = H+Q_1$. So $C$ is contained in a
translate of $H$ in $A$, which is impossible. 

Hence  $\mathrm{Stab}(C-C)$ contains an elliptic curve $E$. Write $f \colon A \rightarrow A/E$. Then $\dim (f(C) - f(C)) = \dim f(C-C) = 1$. Thus the stabilizer of $f(C)$ contains $f(C) - f(C)$ and in particular, $f(C)$ is a coset, which must be of dimension $1$. Hence $C$ is contained in a coset of an abelian surface in $A$. This contradicts $\dim A \ge 3$ and completes the proof. 
\end{proof}

\section{A pencil of curves: N\'{e}ron-Tate distance between algebraic points}\label{SectionNeronTateDistance}

Let $k$ be a number field with algebraic closure $\overline k$ and
suppose $S$ is a smooth,  irreducible,  affine curve defined
over $\overline k$.
We keep the setup from the end of $\mathsection$\ref{sec:vojta}. So 
 $S$ is embedded in some projective space and the absolute logarithmic Weil height pulls back to a height function $h\colon S(\overline k)\rightarrow \IR$.
Moreover, $\pi \colon \cA \rightarrow S$ is an abelian scheme embedded
in a suitable manner
in $\IP^n\times S$.
Let  $\hat h \colon \cA(\overline k)\rightarrow
[0,\infty)$ denote the fiberwise N\'eron--Tate height.
Finally,
 $\cC \subset \cA$ is an irreducible closed surface that dominates
$S$
and such that the   generic fiber of $\pi|_{\cC}\colon \cC \rightarrow S$ is geometrically irreducible.
For $s\in S(\overline k)$ we write $\cA_s = \pi^{-1}(s)$ and $\cC_s = \pi|_{\cC}^{-1}(s)$. 

  



Let $K$ be any field extension of $\overline k$ and suppose $A$ is an
abelian variety defined over $K$. The $K/\overline k$-trace of $A$ is
a final object in the category of pairs $(H,\phi)$ where $H$ is an
abelian variety defined over $\overline k$ and $\phi\colon
H\otimes_{\overline k} K\rightarrow A$ is a homomorphism of abelian
varieties. As $k$ has characteristic $0$ the $K/\overline k$-trace $\tr{K/\overline
  k}{A}$ of $A$ exists and the canonical homomorphism $\tr{K/\overline
  k}{A}\otimes_{\overline k}K\rightarrow A$ is a closed immersion, we
refer to \cite{lang:av} for these and other facts on the trace.
We consider $\tr{K/\overline
  k}{A}\otimes_{\overline k}K$ as a subvariety of $A$.

Let $K$ be the function field $\overline{k}(S)$ of $S$. The generic
 fiber of $\cC \rightarrow S$ is the geometrically irreducible curve $\cC_\eta$ over $K$. The generic fiber $A$ of $\cA\rightarrow S$ is an
abelian variety defined over $K$. We fix an algebraic closure
$\overline K$ of $K$ and write $A_{\overline K} =
A\otimes_K\overline{K}$
and $(\cC_\eta)_{\overline K} = \cC_\eta \otimes_K\overline K$. The goal of this section is to
prove the following proposition.
\begin{proposition}
  \label{lem:ballcount}
   Assume $(\cC_\eta)_{\overline K}$ is not contained in a  coset of 
   $\tr{\overline K/\overline   k}{A_{\overline K}}\otimes_{\overline k}\overline K$
   in    $A_{\overline K}$.
There exist constants $c_2,c_3\ge 1$ that depend only on $\cC$ with the following property. 
If $s \in S(\overline k)$ satisfies $h(s)\ge c_3$ then $\cC_s$ is integral. 
If in addition $\cC_s$ is not a coset in $\cA_s$, then 
  \begin{equation}
    \label{eq:pointsinball}
\# \left\{ Q \in \cC_s(\overline k) : \hat h(Q-P)\le \frac{h(s)}{c_3}\right\}
\le c_2 \quad \text{ for all }\quad P \in \cC_s(\overline k).
  \end{equation}
\end{proposition}


\subsection{Preliminary setup}\label{SubsectionSetupNTDistance}
We keep the notation from above. We begin by making two reduction
steps. Suppose $S'$ is a smooth, irreducible, affine curve that is finite over
$S$ and let $\cA_{S'} = \cA\times_S S'$ and $\cC_{S'} = \cC\times_S S'$. If $s'\in S'(\overline k)$ lies above a
point $s\in S(\overline k)$, then we can identify $(\cA_{S'})_{s'}$ with
$\cA_s$ and $(\cC_{S'})_{s'}$ with $\cC_s$. If $S$ and $S'$ are both
embedded in some, possible different, projective spaces, then by basic height theory we can bound
$h(s)$ from below linearly in terms of $h(s')$. To prove
Proposition~\ref{lem:ballcount} for $\cC$ it suffices to prove it
for $\cC_{S'}$.

First, this observation  allows us to reduce to hypothesis 
\begin{center}
  {\tt (H1):} All endomorphisms of $A_{\overline K}$ are
  endomorphisms of $A$.
\end{center}
Indeed, all geometric endomorphisms of $A$ are defined over a finite field extension $K'/K$.
By a result of Silverberg \cite{Silverberg:92}  and  in characteristic zero it suffices to choose $K'$ such
that all $3$-torsion points of $A$ are $K'$-rational.  Moreover, there
is a smooth, irreducible, affine curve $S'$ and a finite morphism $S'\rightarrow S$ that corresponds to $K'/K$.

It follows in particular that the canonical morphism $\tr{K/\overline
  k}{A}\otimes_{\overline k} \overline K \rightarrow \tr{\overline
  K/\overline k}{A_{\overline K}} \otimes_{\overline
  k}\overline K$ is an isomorphism. So the trace of $A$ does not
increase when replacing $K$ by an algebraic extension of itself.


Second, it suffices to prove the proposition under hypothesis
\begin{center}
{\tt (H2):} 
The curve $\cC_\eta$ is
  not contained in any proper coset of $A$.
\end{center}
Indeed, suppose there is a finite extension $K'/K$ such that
$(\cC_\eta)_{K'} \subset \sigma + A'$ where $\sigma \in A(K')$ and where
$A'$ is an abelian subvariety of $A_{K'}$ of minimal dimension with this property.
Let $S'$ be a smooth, irreducible, affine curve with a finite morphism $S'\rightarrow S$ that corresponds to $K'/K$.
The Zariski closure $\cA'$ of $A'$ in $\cA\times_S S'$ is its N\'eron model and 
the Zariski closure  of $(\cC_{\eta})_{K'}-\sigma$ is a new surface
$\cC'\subset \cA'$. If $\dim A' < \dim A$ we use functorial properties of
the height machine for the N\'eron--Tate height and apply Proposition 
\ref{lem:ballcount} to $\cC'$ by induction on $\dim A$.


For the rest of this section we assume that ${\tt (H1)}$ and ${\tt
  (H2)}$ hold true. 



Let $H = \tr{K/\overline k}{A}$, it is an abelian variety over
$\overline k$. We identify $H_K=H\otimes_{\overline
  k}K$ with the image of the closed immersion $H_K\rightarrow A$, it is an abelian subvariety of
$A$. We fix an abelian subvariety $G$ of $A$ with
$H_K+G=A$ such that $H_K\cap G$ is finite.
The latter condition implies 
 $\tr{K/\overline k}{G}=0$.
Moreover,  addition induces an isogeny $H_K\times G\rightarrow A$. There
is an isogeny $A\rightarrow H_K\times G$ going in the reverse
direction that, when composed with addition, is multiplication by a
non-zero integer on $A$. 
We write $\psi\colon A\rightarrow H_K$ for the
composition of  $A\rightarrow H_K\times G$ followed
by the projection to $H_K$. Then $\psi|_{H_K}\colon H_K\rightarrow
H_K$ is an isogeny and  $\psi(G)=0$.

Note that {(\tt H1)}
holds for $A^2$, in particular, the trace of $A^2$ does not increase
after a finite field extension of $K$.

We write $\psi\times\psi$ for the square
$A^2\rightarrow H_K^2$ of $\psi$; it is the composition of an isogeny
$A^2\rightarrow H_K^2\times G^2$ followed by the projection to $H_K^2$. 
There
is no non-zero homomorphism between $H_K^2$ and $G^2$ as the former comes from an abelian variety of $\overline k$ and
the latter
has $K/\overline k$-trace zero.
So any abelian subvariety  of
$H_K^2 \times G^2$ is a product of an abelian subvariety of $H_K^2$
and an abelian subvariety of $G^2$. We can thus decompose any  abelian
subvariety  of $A^2$ into a sum $B'+B''$ of abelian subvarieties
$B'\subset H_K^2$ and $B''\subset G^2$ such that
$\psi\times\psi(B'')=0$.

The Zariski closure of $H_K^2\subset A^2$ in $\cA\times_S \cA$ is $S\times H^2$
where we consider the abelian variety $H^2$ over $\overline k$ as being contained in each fiber $\cA^2_s$ where
$s\in S(\overline k)$. Our $\psi\times\psi$ from above extends from the generic
fiber to  a homomorphism $\cA \times_S \cA\rightarrow
S\times H^2$ of abelian schemes over $S$, which we still denote by $\psi\times\psi$. On each fiber $\cA_s^2$ it
is the square of a homomorphism $\cA_s \rightarrow H$
whose restriction to $H$ is an isogeny.

Next consider the proper morphism on the family
  \begin{equation*}
    \diffS \colon \cA\times_S \cA \times_S \cA \rightarrow \cA \times_S \cA
  \end{equation*}
defined fiberwise via $\diffS(P_0,P_1,P_2) = (P_1-P_0,P_2-P_0)$ for all
$(P_0,P_1,P_2)\in (\cA\times_S \cA \times_S \cA)(\overline k)$. 
Denote by $\cX = \diffS (\cC\times_S \cC \times_S \cC)$. 
Then $\cX$ is an irreducible, closed subvariety of  $\cA\times_S\cA$.
We will apply the results of $\mathsection$\ref{sec:avoid} to the fibers of $\cX$.


Let $\cX^*$ be as after \cite[Definition~1.2]{GaoHab}.
The structure of $\cX^*$ is clarified by 
\cite[Proposition~1.3]{GaoHab}. In particular, $\cX\ssm \cX^*$ is Zariski closed in
$\cX$. But more is true under {\tt (H1)}, there exist abelian subvarieties
$B_1,\ldots,B_t \subset A^2$  such that 
$\cX \ssm \cX^*$ restricted to the generic fiber is a finite union
$\bigcup_{i=1}^t (Z_i + B_i + \sigma_i)$ where $Z_i$ is an irreducible closed
subvariety of $H^2$
and  $\sigma_1,\ldots,\sigma_t$ are images of torsion points under $A_{\overline K}\rightarrow A$. 

We now demonstrate that we may assume $\psi\times\psi(B_i)=0$ for all $i$.
As above we can decompose
$B_i = B'_i + B''_i$ with $B'_i\subset H_K^2$  and $B''_i \subset
G^2$. Then $Z_i+B'_i\subset H_K^2$ and  $\psi\times\psi(B''_i)=0$.
Note that $Z_i+B_i = (Z_i+B'_i)+B''_i$. 
So after replacing $Z_i$ with $Z_i+B'_i$ and $B_i$ with $B''_i$ we may
assume $\psi\times\psi(B_i)=0$. 

The Zariski closure $\cB_i$ of  $B_i$ in $\cA$ is a subvariety
of $\cA$ for all $i$. Then $\cB_i$ is the N\'eron model of $B_i$. In
particular, the  fibers
$\cB_{i,s}$ are abelian subvarieties of  $\cA_s$ for all $s\in S(\overline k)$. Moreover, we
have $\psi\times\psi(\cB_i) = 0$ after having done the modification
above. 
Let $\cZ_i$ denote the Zariski closure of $Z_i$ in
$\cA$. After possibly increasing $t$ each fiber $\cZ_{i,s}$ is a finite union of irreducible components
$Z_{i,s}$ where $i \in \{1,\ldots, t\}$.

Let us recapitulate. After possibly increasing $t$  we have for all
$s\in S(\overline k)$ that
\begin{equation}
  \label{eq:Xstarcomplement}
 (\cX\ssm \cX^*)_s= \bigcup_{i=1}^{t} \left(Z_{i,s}+\cB_{i,s}+T_{i,s}\right)
\end{equation}
where 
$Z_{i,s}\subset H^2$ is irreducible,
$\cB_{i,s}$ is an abelian subvariety of $\cA_s^2$
with $\psi\times\psi(\cB_{i,s})=0$,
and $T_{i,s}$ is a torsion point on $\cA_s^2$ for all $i$.

\subsection{Proof of Proposition~\ref{lem:ballcount}}
Our main tool to prove Proposition~\ref{lem:ballcount} is a  height lower bound by the second- and
third-named authors \cite{GaoHab} which we state below for the
particular subvariety $\cX$ of $\cA\times_S\cA$. 
It replaces height lower bounds developed by
David--Philippon \cite{DPvarabII,DaPh:07}.
\begin{theorem}
  \label{lem:ghapplication}
  There exists a constant $c_1\ge 1$ depending on $\cC$ but independent
  of $s$, such that if
 $P,P_1,P_2 \in \cC_s(\overline k)$ with
$\diffS(P,P_1,P_2) \not\in \bigcup_{i=1}^{t} \left(Z_{i,s}+\cB_{i,s}+T_{i,s}\right)$, then 
\begin{equation}
\label{eq:thm14applied}
  h(s) \le c_1 \max\{1,\hat h(P_1-P) ,\hat h(P_2-P)\}. 
\end{equation}
\end{theorem}
\begin{proof}
  This is \cite[Theorem~1.4]{GaoHab} applied to $\cX$ together with \eqref{eq:Xstarcomplement}.
\end{proof}


We keep the notation from the previous subsection.

By hypothesis {\tt (H2)}, the curve $\cC_{\eta}$ is
not contained in any proper coset of $A$. Thus the $\dim(A)$-fold
sum of $\cC_{\eta} - \cC_{\eta}$ equals $A$. This sum lies
Zariski dense in $\cA$. Moreover, it contains the $\dim(A)$-fold sum of
$\cC_s-\cC_s$ for all $s$ as
addition $\cA\times_S \cA\rightarrow \cA$ and inversion
$\cA\rightarrow \cA$ are proper morphisms. So the $\dim A$-fold sum of
$\cC_s-\cC_s$ equals $\cA_s$ for all $s$. It follows that $\cC_s$ is
not contained in any proper coset in $\cA_s$ for all $s$.

Suppose $H=\cA_s$ for some $s\in S(\overline k)$. Then $A$ is
constant, \textit{i.e.},  $A = H$. But this contradicts the hypothesis
in Proposition \ref{lem:ballcount} on
$(\cC_\eta)_{\overline K}$. So $H\not\subset \cA_s$ for all $s$.

Note that $\cC_s$ is integral for all but finitely many $s\in
S(\overline k)$ as $\cC_\eta$ is geometrically irreducible. We assume throughout that $h(s)\ge c_3$ with $c_3$ large
enough to ensure that $\cC_s$ is integral. We may assume $c_3>
c_1$ with $c_1$ from Theorem~\ref{lem:ghapplication}

As in Proposition \ref{lem:ballcount} suppose that 
 $\cC_s$ is not a coset in $\cA_s$; in particular
$\cC_s\not=\cA_s$. 
Above we saw that $\cC_s$ is not contained in any proper coset in $\cA_s$.
So for any $P \in \cC_s(\overline k)$, we may apply
Corollary~\ref{LemmaAvoid} to $\cC_s$ and all $Z=Z_{i,s}, B=\cB_{i,s},
(Q_1,Q_2) = T_{i,s}$ with $i\in \{1,\ldots,t\}$. We get $(\cC_s-P)^2 \not\subset \bigcup_{i=1}^{t}
\left(Z_{i,s}+\cB_{i,s}+T_{i,s}\right)$ as $(\cC_s-P)^2$ is
irreducible. 

Let $q \colon \cA \times_S \cA \rightarrow \cA$ be the projection to the first factor. Let $Y$ be an irreducible component of $(\cC_s - P)^2 \cap \bigcup_{i=1}^t(Z_{i,s}+\cB_{i,s}+T_{i,s})$. Then $\dim Y\le 1$.
As $t$ is bounded uniformly in $s$ and by B\'ezout's Theorem the number of
$Y$ 
is bounded by $k_1 \ge 1$, which is independent of $s$.
If $\dim q(Y) = 1$, then each fiber of $q|_Y$ has dimension $0$.
The degrees of $\cC_s,Z_{i,s},$ and $\cB_{i,s}$ are uniformly bounded
in $s$.
So again by B\'ezout's Theorem  we find that the cardinality of each fiber
of $q|_Y$ is bounded  from above, say by $k_2 \ge 1$, independently of $s$. 
Let $c_2 = k_1k_2$.



Now assume that the cardinality of the set displayed in
(\ref{eq:pointsinball}) is strictly larger than $c_2$.
As $c_2\ge k_1$ we max fix $P_1$ from this set such that $P_1-P$ is not
equal to a zero dimensional $q(Y)$.
Then as $c_2\ge k_1k_2$ we can find $P_2$ from the said set
such that $(P_1-P,P_2-P)$ does not lie in any $Y$ as above. 
Thus
$(P_1-P,P_2-P) = \diffS(P,P_1,P_2) \not\in \bigcup_{i=1}^{t} \left(Z_{i,s}+\cB_{i,s}+T_{i,s}\right)$. So Theorem~\ref{lem:ghapplication} implies $h(s) \le c_1 \max\{1,\hat h(P_1-P),\hat h(P_2-P)\}$. Recall that $h(s)\ge c_3> c_1>0$. So $h(s) \le c_1 h(s)/c_3$; here we used that $P_1$ and $P_2$ are points in the set displayed in (\ref{eq:pointsinball}). This is a contradiction. \qed

\section{A pencil of curves: the desired bound}


Let $\pi \colon \cA \rightarrow S$, $\cC \subset \cA$, $h \colon S(\overline k) \rightarrow \IR$ and $\hat h \colon \cA(\overline k)\rightarrow [0,\infty)$ be as in $\mathsection$\ref{SectionNeronTateDistance}.
  Note that we still assume that $S$ is affine.
  

The goal of this section is to prove the following version of
Theorem~\ref{thm:main}.

\begin{theorem}
\label{thm:familymain}
There exists a constant $c=c(\pi,h,\cC)\ge 1$ with the following
property. Let $s\in S(\overline k)$ such that $h(s) \ge c$. Then 
 $\cC_s$ is integral. If in addition $\cC_s$ is not a coset of
$\cA_s$ and  if $\Gamma$ is a finite rank subgroup of $\cA_s(\overline k)$
of rank $\rho$, then $\# \cC_s(\overline k)\cap\Gamma \le c^{1+\rho}$.
\end{theorem}

As in $\mathsection$\ref{SectionNeronTateDistance} we let $K=\overline k(S)$ and fix an algebraic closure $\overline K$ of $K$. Again, $A$ is the generic fiber of $\cA\rightarrow S$ and $\cC_\eta$ is the generic fiber of $\cC\rightarrow S$.


\subsection{The non-isotrivial case}



Let $s \in S(\overline k)$, and let $\Gamma$ be a finite rank subgroup of $\cA_s(\overline k)$ of rank $\rho$.
\begin{lemma}
\label{lem:largepoints}  
There exists a constant $c_1\ge 1$ depend on $\cC$ but independent of $s$
with the following property. 
If $\cC_s$ is not a coset in $\cA_s$, then  $\#\{ P \in \cC_s(\overline k) \cap\Gamma
: \hat h (P) > c_1\max\{1,h(s) \}\} \le c_1^{\rho}$.
\end{lemma}
\begin{proof}
  This is just Corollary~\ref{cor:countlargepts2}.
\end{proof}

We will apply the following packing lemma where the norm
 $|\cdot|$  on  the $\rho$-dimensional  $\IR$-vector space $\Gamma\otimes_\IZ\IR$ is
 induced by $\hat h^{1/2}$. 

\begin{lemma}
\label{lem:spherecovering}  
  Let $0<r\le R$ 
and   let $M\subset \Gamma\otimes_\IZ\IR$ be a subset of the ball of radius $R$ around $0$.
  There exists a finite set $\Sigma \subset M$
  with $\#\Sigma \le (1+2R/r)^\rho$ such that any point of $M$
  is contained in  a closed ball of radius $r$ center at some element of $\Sigma$.
  \end{lemma}
\begin{proof}
This follows from R\'emond \cite[Lemme~6.1]{Remond:Decompte}. 
\end{proof}

Now we are ready to prove Theorem~\ref{thm:familymain}
if
$(\cC_\eta)_{\overline K}$ is not contained in any coset of
   $\tr{\overline K/\overline   k}{A_{\overline K}}\otimes_{\overline k}\overline K$
   in    $A_{\overline K}$.
So the hypothesis of Proposition~\ref{lem:ballcount} is fullfilled.

Let $c_1$ be as in Lemma~\ref{lem:largepoints}, and let $c_2,c_3$ be as in Proposition~\ref{lem:ballcount}; they are independent
of $s$.
By Lemma~\ref{lem:largepoints} there are at most $c_1^\rho$ points in
$\cC_s(\overline k)\cap\Gamma$ of N\'eron--Tate height  strictly greater than $c_1\max\{1,h(s)\}$.
Our goal is to get a similar bound for the points in $\cC_s(\overline k)\cap\Gamma$
of N\'eron--Tate height at most $c_1\max\{1,h(s)\}$.
In doing this we may assume that $h(s)\ge c=c_3\ge 1$. 

To count the remaining points we combine Proposition~\ref{lem:ballcount}
and Lemma~\ref{lem:spherecovering}. More precisely we apply Lemma~\ref{lem:spherecovering}  to
$R
= (c_1h(s))^{1/2}$ and $r = \left(h(s)/c_3\right)^{1/2}$.
The dependency on $h(s)$ cancels out in $R/r
= (c_1 c_3)^{1/2}\ge 1$.
For $M$ we take the image in $\Gamma\otimes_\IZ\IR$ of all $P\in \cC_s(\overline k)\cap\Gamma$ such
that $\hat{h}(P) \le R^2$.  We can cover the closed ball in $\Gamma\otimes_\IZ\IR$ of radius $R$
around $0$ using at most $(1+2(c_1c_3)^{1/2})^\rho$ closed balls of radius
$r$ centered around certain points coming from $\cC_s(\overline k)\cap\Gamma$. By
Proposition~\ref{lem:ballcount} at most $c_2$ points of $\cC_s(\overline k)\cap \Gamma$ end
up in one of the closed balls of radius $r$.
So the number of points in $\cC_s(\overline k)\cap\Gamma$ of height at most
$c_1h(s)$ is bounded from above by
\begin{equation*}
c_2 (1+2(c_1c_3)^{1/2})^{\rho}. 
\end{equation*}
Adding the cardinality bounds for the large and small points
yields $\# \cC_s(\overline k)\cap\Gamma \le c_1^{\rho}+  c_2(1+2(c_1c_3)^{1/2})^{\rho}$.
The theorem follows as $c_1,c_2,$ and $c_3$ are independent of $s$.

\subsection{The isotrivial case}
Suppose that 
$(\cC_\eta)_{\overline K}$ is contained in a coset of
   $\tr{\overline K/\overline   k}{A_{\overline K}}\otimes_{\overline k}\overline K$
in    $A_{\overline K}$.
Roughly speaking, all fibers $\cC_s$ are contained in the same abelian
variety and we can refer to the constant case. 
We use the notation in $\mathsection$\ref{SubsectionSetupNTDistance}.
With similar argument below Theorem~\ref{lem:ghapplication}, we may
and do reduce to the case {\tt (H1)}. 
In particular, $\tr{\overline K/\overline   k}{A_{\overline K}}\otimes_{\overline k}\overline K = 
    A_{\overline K}$.
We may assume that  $\cA = S \times H$. 
Note that each fiber $\cC_s \subset H$ has uniformly bounded degree as $s$ varies over $S(\overline k)$. 
Then the desired bound follows from R\'emond
\cite[Th\'eor\`eme~1.2]{Remond:Decompte} as the ambient abelian
variety $H$ is independent of $s$. This completes our proof of Theorem~\ref{thm:familymain}.   \qed

\subsection{Proof of Theorem~\ref{thm:main}}
We may freely remove finitely many points from $S$. In particular we
may assume that $S$ is affine.
Let $\cC_\eta$ denote the generic fiber of $\cC\rightarrow S$.
There
is a smooth, geometrically irreducible curve $S'$ and a finite
morphism $S'\rightarrow S$ such that $\cC_\eta$ has a $\overline
k(S')$-rational point. So $S'$ is affine too. We use this point to embed $\cC_\eta
\otimes_{\overline k(S)} \overline k(S')$ into its Jacobian. Let $\cC'
= \cC\times_S S'$ and let $\cJ$ denote the relative Jacobian of
$\cC'$. By the N\'eron mapping property we get a closed immersion
$\cC'\rightarrow \cJ$ over $S'$, after again possibly shrinking $S$.
If $s'\in S'(\overline k)$ lies above $s\in S(\overline k)$, we can
identify the fiber $\cC'_{s'}$ with the fiber $\cC_s$. 
Note that $\cC'_{s'}$ is smooth and projective
of  genus $g \ge 2$, 
  so it cannot be a coset in the fiber $\cJ_{s'}$. 
We obtain the desired bound in  Theorem~\ref{thm:main}  from
Theorem~\ref{thm:familymain} for an embedding of $\cC_s$ into
$\jac{\cC_s}$ using \emph{some} base point. But the cardinality bound
actually holds for any base point as we can enlarge $\Gamma$ by
adding an additional generator and replacing $c$ by $c^2$. 
\qed

\section{An example involving hyperelliptic curves}
This section contains   examples for the reader who does not wish to
disappear completely in the realm of abelian
varieties.

Let us consider a squarefree polynomial $Q\in \IQ[x]$ of
degree $d-1\ge 4$. We  will bound   the number of rational solutions of the  family
\begin{equation*}
 y^2 = (x-s) Q(x)
\end{equation*}
in terms of the rational parameter
$s$ with $Q(s)\not=0$.

The equation determines a one-parameter family of hyperelliptic curves
$\cC$.
The genus $g$ of each member is $(d-2)/2$ for $d$ even and $(d-1)/2$ for $d$ odd. 
For  rational $s$ with $Q(s)\not=0$,
the point $(s,0)$ is a  rational Weierstrass point of $\cC_s$.

We can bound the rank $\jac{\cC_s}(\IQ)$ from above using
\cite[Theorem~C.1.9]{DG2000}. To apply this estimate we  pass to an extension $k/\IQ$
over which all $2$-torsion points of $\jac{\cC_s}$ are defined. The
difference of Weierstrass points of $\cC_s$ generate the $2$-torsion in the
Jacobian. If $d$ is even, the Weierstrass points come from roots of
$(x-s)Q(s)$. If $d$ is odd, there is an additional rational
Weierstrass point at infinite. We take for $k$ the splitting field
of $Q$.

For $s\in \IQ$, let $\Delta(s)$ denote the discriminant of $(x-s)Q(x)
\in \IQ[x]$. Then $\Delta(s)$ is a polynomial in $s$ with rational
coefficients and degree at most $2d-2$.
Suppose now $Q(s)\not=0$. 
Then  $\Delta(s)\not=0$ as $Q$ is squarefree.
Let $q\in\IN=\{1,2,3,\ldots\}$ such that $qQ \in \IZ[x]$.
Say  $s=a/b$ with $a,b$ coprime integers and $b\ge 1$, then
\begin{equation}
  \label{eq:discgeneral}
    (bq)^{2d-2} \Delta(s) \in\IZ. 
\end{equation}
Let $\cP$ be the set of primes $p$ with $p\mid 2bq$ or
such that $p$ divides (\ref{eq:discgeneral}). 
Any prime number outside of $\cP$  is a prime where $\cC_s$ and its Jacobian have
good reduction. 
We estimate $\#\cP \le \omega(2bq) + \omega((bq)^{2d-2} \Delta(s))$ where $\omega(n)$
denotes the number of prime numbers dividing $n$.
Observe that $|b^{2d-2}\Delta(s)|=O(\max\{|a|,b\}^{2d-2})$ where
here and below the implicit constant is allowed to depend on $Q$ and
$k$ but not on $s$. 
Using $\omega(n) = O(\log n/\log\log n)$ for $n\ge 3$  we obtain
\begin{equation}
  \label{eq:calPbound}
  \#\cP =O( \log H^*(s) / \log\log H^*(s)) 
\end{equation} 
 where $H^*(s) = \max\{3,e^{h(s)}\}=\max\{3,|a|,b\}$.

We take for $\cP'$ the set of maximal ideals in the ring of
integers of $k$ that contain a prime in $\cP$.
So $\#\cP' \le [k:\IQ] \#\cP$. 
Next we add finitely many maximal ideals to $\cP'$ such that the
 ring of $\cP'$-integers in $k$ is a principal ideal domain.
 The number of ideals to add can be bounded from above solely in terms of $k$. 
Recall that $\jac{\cC_s}(k)$
contains the full $2$-torsion subgroup of $\jac{\cC_s}$ (of order $2^{2g}$).
By \cite[Theorem~C.1.9]{DG2000} with $m=2$ and $r_1+r_2\le
[k:\IQ]$ the rank $ \rk{\jac{\cC_s}(k)}$ is at most
\begin{equation}
  \label{eq:rkjacbound}
  2 g (  [k:\IQ]-1+\# \cP') = O( \log H^*(s) / \log\log H^*(s)).
\end{equation}
having used (\ref{eq:calPbound}).


 \begin{corollary}
   \label{cor:1234s}
   Let $Q\in \IQ[x]$ be as above. 
There exists a constant  $c=c(Q)\ge 1$ such that
\begin{equation*}
  \# \left\{ (x,y)\in\IQ^2 : y^2 = (x-s) Q(x) \right\}
  \le  H^*(s)^{c/\log\log H^*(s)}
\end{equation*}
for all $s\in \IQ$ with $Q(s)\not=0$. 
\end{corollary}
\begin{proof}
  This is a direct application of Corollary~\ref{cor:main}
  applied to
 the Mordell--Weil group
$\Gamma = \jac{\cC_s}(\IQ)$ 
  whose rank is bounded by
  (\ref{eq:rkjacbound}). The base $S$ in the corollary is the affine
  line punctured at the roots of $Q$. 
\end{proof}

We conclude that the number of rational points on $\cC_s$ grows subpolynomially in the
\emph{exponential} height of the parameter $s$.

Let us conclude with a few comments regarding the special case where
$Q = x(x-2)(x-6)(x-8)(x-12)(x-20)$; the results below extend without
much effort to more general $Q$.  The hyperelliptic curve is presented by
\begin{equation}
  \label{eq:0268}
  y^2 = (x-s)x(x-2)(x-6)(x-8)(x-12)(x-20),
\end{equation}
and 
$\Delta(s)= 2^{48} 3^{14}5^4 7^2 s^2(s-2)^2(s-6)^2(s-8)^2(s-12)(s-20)^2$.
Here the genus is $g=3$ and there is a single rational point at
infinity. We can take $k=\IQ$ as
$Q$ splits completely over the rationals.

Given $s=a\in \IZ\ssm\{0,2,6,8,12,20\}$ we can take for
$\cP$  the primes dividing
$\Delta(s)$; this includes $2$. 
Here the choice, $\cP'=\cP$ is possible.
The Mordell--Weil rank of the $\IQ$-rational points of the Jacobian
is at most $2g \#P \le 6\omega(\Delta(s))$ by
\cite[Theorem~C.1.9]{DG2000}. Thus the number of $(x,y)\in
\IQ^2$ satisfying (\ref{eq:0268}), for fixed $s\in\IZ\ssm\{0,2,6,8,12,20\}$, is at most
$c^{\omega(s(s-2)(s-6)(s-8)(s-12)(s-20))}$
where $c\ge 1$ is a constant independent of $s$.

For all primes $p$, the set $\{0,2,6,8,12,20\}$ does not map
surjectively to $\IZ/p\IZ$. By a result of Halberstam and Richert
\cite[Theorem 4]{HR:distpolyseq} there exists $t\ge 1$ and infinitely
many integers $s$ with $\omega(s(s-2)(s-6)(s-8)(s-12)(s-20))\le t$.
Our result implies that the number of rational solutions of
(\ref{eq:0268}) is uniformly bounded for infinitely many $s\in\IZ$.
Generalizing to a different set of roots requires an affine linear
change of coordinates. 

The normal order of $\omega(s)$ is $\log\log s$ by the
Hardy--Ramanujan Theorem. This theorem and our result imply that there
exist $A\ge 1$ and $N\subset\IN$ of natural density $1$ such that the
number of $(x,y)\in\IQ^2$ satisfying (\ref{eq:0268}) is at most $(\log
s)^A$ for all $s\in N$.
\bibliographystyle{amsplain}

\vfill

\end{document}